\documentclass{amsart}
\usepackage{amssymb,amsmath}
\usepackage{enumitem}
\usepackage{xcolor}
\usepackage{array}
\usepackage{longtable}
\usepackage{tikz}
\usepackage{hyperref}
\hypersetup{hidelinks, pdftitle={Improving Razmyslov's (1974) bound: the Nagata-Higman index is at most n\textasciicircum 2-2 for n>=3}, pdfauthor={Dave Goel}}

\newtheorem{theorem}{Theorem}[section]
\newtheorem{proposition}[theorem]{Proposition}
\newtheorem{lemma}[theorem]{Lemma}
\newtheorem{corollary}[theorem]{Corollary}

\newtheorem{mainthm}{Theorem}

\theoremstyle{definition}

\theoremstyle{remark}
\newtheorem{remark}[theorem]{Remark}
\newtheorem*{example}{Example}

\newcommand{\sgn}{\operatorname{sgn}}
\newcommand{\End}{\operatorname{End}}
\newcommand{\Res}{\operatorname{Res}}
\newcommand{\Stab}{\operatorname{Stab}}
\newcommand{\Sym}{\operatorname{Sym}}
\newcommand{\tr}{\operatorname{tr}}
\newcommand{\Arc}{\operatorname{Arc}}
\makeatletter
\def\@cite#1#2{{\@citestyle[\citeform{#1}\if@tempswa,~#2\fi]}}%
\makeatother

\title[Improving Razmyslov's bound]{Improving Razmyslov's (1974) bound:\\ the Nagata--Higman index is at most $n^2-2$ for $n\ge3$}
\author{Dave Goel}
\email{david44hello@gmail.com}
\thanks{The Lean~4 proof of this paper's results:\newline
\url{https://github.com/david44hello/papers_supporting_materials/tree/nagata-higman-index-v1.12/nagata-higman-index}.}
\subjclass[2020]{Primary 16R10; Secondary 16N40, 16R30, 20C30}
\keywords{Nagata--Higman theorem, nil algebras, nilpotency index, Kuzmin's conjecture, polynomial identities, symmetric group}
\date{November 1, 2026}

\begin{document}

\begin{abstract}
Let $d(n)$ be the least $m$ such that every associative (not necessarily unital) algebra over a field of characteristic zero
satisfying the identity
$x^n=0$ also satisfies $x_1x_2\cdots x_m=0$; this is the Nagata--Higman index. Razmyslov proved $d(n)\le n^2$ in 1974; Kuzmin
proved $d(n)\ge n(n+1)/2$ and conjectured equality; this is the Kuzmin (1975) conjecture. We prove $d(n)\le n^2-2$ for $n\ge3$ (and, by
a simpler case of the same method, $d(n)\le n^2-1$ for $n\ge2$); in particular $d(6)\le34$. To our knowledge, $n=6$ is the first case in
which the new bound improves the known bounds, and this is the first improvement of Razmyslov's bound for general $n$. For $m=n^2-1$ and
$m=n^2-2$, Formanek's criterion reduces $d(n)\le m$ to the statement that no nonzero linear combination of $(m+1)$-cycles lies
in one particular isotypic block of the group algebra of $S_{m+1}$. For $m=n^2-1$ this statement follows from the
branching rule; for $m=n^2-2$ we prove it using the Jucys--Murphy elements and a theorem of Frobenius on the characters of
the alternating group. These bounds have also been
proved in the Lean proof assistant.
\end{abstract}

\maketitle

\section{Introduction}\label{sec:intro}

Throughout, $F$ is a field of characteristic zero. An associative (not necessarily unital) $F$-algebra is \emph{nil of index
$n$} if $x^n=0$ for each of its elements $x$.
By the Nagata--Higman theorem, such an algebra is nilpotent: there is an $m$, depending only on $n$, such that
$x_1x_2\cdots x_m=0$ for all elements $x_1,\dots,x_m$ of every such algebra. The least such $m$ is the \emph{Nagata--Higman index}
$d(n)$. Razmyslov \cite{Ra74} proved $d(n)\le n^2$ in 1974; Kuzmin \cite{Ku75} proved $d(n)\ge n(n+1)/2$ in 1975 and conjectured
that equality holds. For $n\ge3$ we lower Razmyslov's bound by two. The table and the list below give our main results, each with the best
result known before this paper.

\begin{table}[t]
\small
\renewcommand{\arraystretch}{1.1}
\begin{center}
\begin{tabular}{|>{\raggedright\arraybackslash}p{0.46\textwidth}|>{\raggedright\arraybackslash}p{0.43\textwidth}|}
\hline
\textbf{New result} & \textbf{Best result before this paper}\\
\hline
$d(n)\le n^2-2$ for $n\ge3$, new for $n\ge6$ (Theorem~\ref{thm:B}) &
$d(n)\le n^2$ (Razmyslov 1974 \cite{Ra74}), still called the best known upper bound in 2022 \cite{Ro22}\\
\hline
$d(6)\le34$ &
$d(6)\le36$ \cite{Ra74}\\
\hline
invariants of $n\times n$ matrices: generated in degree \mbox{$\le n^2-2$} for $n\ge3$, new for $n\ge6$ (Corollary~\ref{cor:invariants}) &
degree $\le n^2$ (\cite{Ra74} with \cite[Thm.~3.3]{Pr76}), ``at the moment the best result'' \cite[Thm.~5.14]{Pr26}\\
\hline
\end{tabular}
\end{center}
\end{table}

\begin{itemize}[leftmargin=*]
\item \emph{Theorem~\ref{thm:A}: $d(n)\le n^2-1$ for every $n\ge2$.} This lowers Razmyslov's bound $n^2$ \cite{Ra74} by one;
Theorem~\ref{thm:B} supersedes it for $n\ge3$. The proof (Section~\ref{sec:A}) uses the square partition $(n,\dots,n)$.
\item \emph{Theorem~\ref{thm:B}: $d(n)\le n^2-2$ for every $n\ge3$.} This lowers Razmyslov's bound \cite{Ra74} by two; with
Kuzmin's lower bound \cite{Ku75}, $n(n+1)/2\le d(n)\le n^2-2$. The proof (Section~\ref{sec:B}) uses the partition \mbox{$(n,\dots,n,n-1)$} and a theorem of Frobenius.
\item \emph{$d(6)\le34$.} The values $d(n)$ for $n\le4$ are known \cite{Du35,Hi56,VL93}, and Higman's recursion with $d(4)=10$
gives $d(5)\le21$ \cite{Hi56}; so $n=6$ is the first case in which Theorem~\ref{thm:B} improves the known bounds.
\item \emph{Corollary~\ref{cor:invariants}: for $n\ge3$, the invariants of $n\times n$ matrices are generated in degree at most
$n^2-2$; this is new for $n\ge6$.} Before, the known bound for general $n$, uniform in the number of matrices, was $n^2$, from Razmyslov's bound and the theorem of Procesi
and Formanek \cite[Thm.~3.3]{Pr76}, \cite[Thm.~6]{Fo86}; Procesi \cite[Thm.~5.14]{Pr26} still calls it ``at the moment the best result''.
\end{itemize}

\subsection{The Nagata--Higman theorem}\label{subsec:NH}
The theorem is due to Dubnov and Ivanov \cite{DI43}, Nagata \cite{Na52} and Higman \cite{Hi56}, and is also called the
Dubnov--Ivanov--Nagata--Higman theorem \cite[\S5]{Dr95}; Higman's proof gives $d(n)\le2^n-1$. The index $d(n)$ is also the least $m$ such that
$x_1\cdots x_m$ lies in the T-ideal generated by $x^n$, that is, in the smallest ideal of the free non-unital algebra
$F\langle x_1,x_2,\dots\rangle$ that contains $x^n$ and is stable under all endomorphisms (here $x$ is one of the letters);
see \cite[\S1.3]{Ro22}. Indeed, the polynomials that vanish under every substitution of elements of an algebra $R$ form a
T-ideal, and it contains $x^n$ if $R$ is nil of index $n$. Conversely, the quotient of $F\langle x_1,x_2,\dots\rangle$ by the
T-ideal of $x^n$ is itself nil of index $n$: each of its elements is the image of a polynomial $h$ without constant term, and
$h^n$ is the image of $x^n$ under the endomorphism that sends $x$ to $h$ and fixes the other letters. In words, a product
vanishes in every algebra that is nil of index $n$ exactly when it is a consequence of the identity $x^n=0$. We use this form in
Section~\ref{sec:criterion}.

\subsection{Bounds and known values}\label{subsec:bounds}
Razmyslov's proof \cite[final remark]{Ra74} uses the trace identities of the matrix algebra $M_n(F)$ and the
Cayley--Hamilton theorem. By Procesi \cite[Thm.~3.3]{Pr76} and Formanek \cite[Thm.~6]{Fo86} (see also \cite{Fo87,Fo91}), the
invariants of any number of $n\times n$ matrices under simultaneous conjugation are generated by traces of products of degree at
most $d(n)$, and this degree bound is sharp when the number of matrices is large enough \cite[Thm.~5.6]{Dr95}. Kuzmin's lower bound \cite{Ku75} comes from the word $x_2x_1x_2x_1^2\cdots x_2x_1^{n-1}$: it has degree \mbox{$n(n+1)/2-1$} and does not lie in the T-ideal of $x^n$
(see \cite[Thm.~1.7]{Do18}), so neither does $x_1x_2\cdots x_{n(n+1)/2-1}$, which maps onto it under an endomorphism.

The exact values $d(1)=1$, $d(2)=3$ and $d(3)=6$ are classical (Dubnov \cite{Du35}; see \cite[\S5]{Dr95} and, for $n=3$,
\cite{Hi56}), and Vaughan-Lee \cite{VL93} computed $d(4)=10$ by computer.
Higman's recursion $d(n)\le2d(n-1)+1$ \cite{Hi56} then gives $d(5)\le21$, and Kuzmin's bound gives $d(5)\ge15$.

To our knowledge, in characteristic zero the bound $n^2$ had not been improved for general $n$ since 1974; see
\cite[\S5]{Dr95}, \cite[\S1.3]{Ro22} and \cite[Thm.~5.14]{Pr26}.
\subsection{Results}
\begin{mainthm}\label{thm:A}
For every $n\ge2$, $d(n)\le n^2-1$.
\end{mainthm}

\begin{mainthm}\label{thm:B}
For every $n\ge3$, $d(n)\le n^2-2$.
\end{mainthm}

For $n=2$, Theorem~\ref{thm:A} only recovers the classical bound $d(2)\le3$. Theorem~\ref{thm:A} is the first, simpler
case of the method: given Formanek's criterion (Proposition~\ref{prop:criterion}), its proof is short
(Lemma~\ref{lem:onecorner}), but we have not found it in print. The gain over Razmyslov's bound is an additive constant; the
gap to Kuzmin's lower bound is still of order $n^2/2$. By the
theorem of Procesi and Formanek (Subsection~\ref{subsec:bounds}), for $n\ge6$ Theorem~\ref{thm:B} also lowers the known degree bound, uniform in the number of matrices, for the
generators of the matrix invariants:

\begin{corollary}\label{cor:invariants}
For $n\ge3$ and every $k\ge1$, the algebra of polynomial invariants of $k$-tuples $(X_1,\dots,X_k)$ of $n\times n$ matrices over
$F$ under simultaneous conjugation is generated by the traces $\tr(X_{i_1}\cdots X_{i_r})$ with $1\le r\le n^2-2$.
\end{corollary}

\begin{proof}
First let $F$ be algebraically closed. By the first fundamental theorem of matrix invariants \cite[Thm.~5.3]{Dr95}, the
invariants are generated by the traces $\tr(X_{i_1}\cdots X_{i_r})$; by the theorem of Procesi and Formanek
(Subsection~\ref{subsec:bounds}), those with $r\le d(n)$ suffice; and $d(n)\le n^2-2$ by Theorem~\ref{thm:B}.
Now let $F$ be arbitrary (of characteristic zero), with algebraic closure $\bar F$. Conjugation preserves degrees, so it
suffices to treat invariants $g$ that are homogeneous of some degree $\delta$ in the $kn^2$ matrix entries. Writing
$\gamma^{-1}=\det(\gamma)^{-1}\operatorname{adj}(\gamma)$ for an invertible matrix $\gamma$, invariance says
$g(\gamma X_1\operatorname{adj}(\gamma),\dots,\gamma X_k\operatorname{adj}(\gamma))=\det(\gamma)^\delta g(X_1,\dots,X_k)$ for
all invertible $\gamma$. Since $F$ is infinite, this holds for all invertible $\gamma$ over $F$ if and only if it holds as a
polynomial identity in the entries of $\gamma$ and of the $X_i$ (the product of $\det(\gamma)$ and the difference of the two
sides vanishes at every point, so it is the zero polynomial, and hence so is the difference); that is a system of linear equations, with integer coefficients, for the coefficients of $g$. So the
invariants of degree $\delta$ over $F$ and over $\bar F$ are the solutions of the same linear system, and the two spaces have
the same dimension. The products of total degree $\delta$ of traces $\tr(X_{i_1}\cdots X_{i_r})$ with $r\le n^2-2$ have integer coefficients and span
the invariants of degree $\delta$ over $\bar F$; so their span over $F$ has the same dimension, and it is the whole space of
invariants of degree $\delta$ over $F$.
\end{proof}

\subsection{Outline of the proofs}\label{subsec:outline}
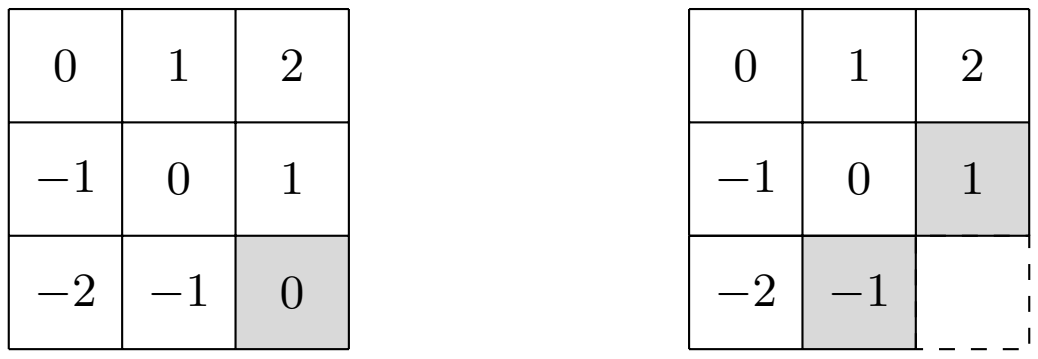
\begin{figure}[t]
\centering
\begin{tikzpicture}[x=0.75cm,y=0.75cm,font=\small]
\fill[black!15] (2,-3) rectangle (3,-2);
\draw[step=1] (0,-3) grid (3,0);
\foreach \i in {1,2,3} \foreach \j in {1,2,3} {\pgfmathtruncatemacro{\c}{\j-\i}\node at (\j-0.5,0.5-\i) {$\c$};}
\node at (1.5,-3.55) {$(3^3)\vdash 9$};
\begin{scope}[xshift=4.5cm]
\fill[black!15] (2,-2) rectangle (3,-1);
\fill[black!15] (1,-3) rectangle (2,-2);
\draw[step=1] (0,-2) grid (3,0);
\draw[step=1] (0,-3) grid (2,-2);
\draw[dashed] (2,-3) rectangle (3,-2);
\foreach \i in {1,2} \foreach \j in {1,2,3} {\pgfmathtruncatemacro{\c}{\j-\i}\node at (\j-0.5,0.5-\i) {$\c$};}
\foreach \j in {1,2} {\pgfmathtruncatemacro{\c}{\j-3}\node at (\j-0.5,-2.5) {$\c$};}
\node at (1.5,-3.55) {$\nu=(3^2,2)\vdash 8$};
\end{scope}
\end{tikzpicture}
\caption{The partitions that fit in the $n\times n$ box in the cases of Theorems~\ref{thm:A} and~\ref{thm:B}, drawn for $n=3$.
Each cell $(i,j)$ shows its content $j-i$; the removable corners are shaded. Left: the square $(3^3)$, the case $N=n^2$; its
only removable corner is the cell $(3,3)$. Right: $\nu=(3^2,2)$, the case $N=n^2-1$, which is the square with the cell $(3,3)$
(dashed) removed; its two removable corners, $(2,3)$ and $(3,2)$, have contents $+1$ and $-1$.}\label{fig:shapes}
\end{figure}

Let $m\ge n$ and $N=m+1$. Let $C_N$ be the span of the $N$-cycles in the group algebra $F[S_N]$. For a partition
$\lambda\vdash N$, let $I_\lambda$ be the isotypic block of $\lambda$, that is, the minimal two-sided ideal of $F[S_N]$ that
acts nontrivially on the irreducible $F[S_N]$-module $S^\lambda$ indexed by $\lambda$ (see (F2)); it acts as zero on the other irreducible modules. We identify $f\in F[S_N]$
with its coefficient function $\sigma\mapsto f(\sigma)$, and we say that $f$ \emph{vanishes on} a set of permutations if
$f(\sigma)=0$ for every $\sigma$ in the set. Our starting point is Formanek's criterion \cite[Thm.~5]{Fo90}; related
statements are in \cite[final remark]{Ra74}, \cite[Thm.~3.3]{Pr76} and \cite{St89}. We prove it in the following cycle form
(Proposition~\ref{prop:criterion}): $d(n)\le m$ if and only if $C_N\cap B_n=0$, where $B_n$ is the sum of the blocks
$I_\lambda$ over the partitions $\lambda\vdash N$ that fit in the $n\times n$ box, that is, have at most $n$ rows and at most
$n$ columns. For example, for $n=2$ only $(2,1)$ fits when $N=3$. This is the nil analogue of the description of the
multilinear trace identities of matrix algebras due to Helling, Razmyslov and Procesi \cite{He74,Ra74,Pr76}; see also
\cite[Thm.~5.7]{Dr95}. The box has $n^2$ cells, so for $N>n^2$ no partition fits, and the criterion gives Razmyslov's bound
(Remark~\ref{rem:razmyslov}). Theorems~\ref{thm:A} and~\ref{thm:B} treat the next two cases, $N=n^2$ and $N=n^2-1$; in each,
exactly one partition fits (Figure~\ref{fig:shapes}).
\begin{enumerate}[label=(\roman*),leftmargin=*, itemsep=1pt, topsep=2pt]
\item The case $N=n^2$ (Theorem~\ref{thm:A}). Only the square $(n^n)=(n,\dots,n)$ fits. An element of $C_N$ vanishes on the
subgroup $S_{N-1}$ of the permutations that fix $N$, since no $N$-cycle fixes a point. The square has only one removable
corner (a cell whose removal leaves the diagram of a partition), so, by the branching rule, the restriction of $S^{(n^n)}$ to
$S_{N-1}$ is irreducible. Hence an element of the block $I_{(n^n)}$ is determined by its coefficients on $S_{N-1}$, and an
element of $C_N\cap I_{(n^n)}$ is zero (Lemma~\ref{lem:onecorner}).
\item The case $N=n^2-1$ (Theorem~\ref{thm:B}). Only $\nu=(n^{n-1},n-1)$, the square with the cell $(n,n)$ removed, fits.
Its two removable corners, $(n-1,n)$ and $(n,n-1)$, have contents $+1$ and $-1$ (the content of the cell $(i,j)$ is $j-i$).
An element $f$ of $C_N$ vanishes on every permutation with a fixed point. Let $\rho_\nu$ be the representation of $F[S_N]$ on
$S^\nu$, and let the \emph{star element} $X^{(j)}:=\sum_{i\ne j}(i\,j)$ be the sum of the $N-1$ transpositions that move
$j$; thus $X^{(N)}=J_N$, where $J_k:=\sum_{i<k}(i\,k)$ is the $k$-th Jucys--Murphy element \cite{Ju74,Mu81,OV96}. The proof has three steps.

First, for $f\in I_\nu$, vanishing on every permutation with a fixed point means that $\rho_\nu(f)$ anticommutes with every
$\rho_\nu(X^{(j)})$ (Lemma~\ref{lem:anticommute}). Indeed, restricted to $S_{N-1}$, the module $S^\nu$ is the sum of two
irreducible modules, one for each corner, on which $\rho_\nu(J_N)$ acts by the contents $+1$ and $-1$. Vanishing on
$S_{N-1}$ says that $\rho_\nu(f)$ has zero diagonal blocks for this sum, that is, that it anticommutes with $\rho_\nu(J_N)$;
conjugation gives the other $j$.

Second, only the scalars commute with every $\rho_\nu(X^{(j)})$ (Lemma~\ref{lem:claimS}): the contents $\pm1$ let us recover
$\rho_\nu(J_{N-1})$, and then $\rho_\nu((N{-}1\;N))$, from the $\rho_\nu(X^{(j)})$, so the $\rho_\nu(X^{(j)})$ generate $\End(S^\nu)$ as an algebra.

Third, since $\nu$ equals its conjugate (transposed) partition, there is an operator $\Gamma$ on $S^\nu$ (over $\mathbb C$),
with $\Gamma^2=\mathrm{id}$, that anticommutes with $\rho_\nu(t)$ for every transposition $t$, and hence with every
$\rho_\nu(X^{(j)})$. So $\Gamma\rho_\nu(f)$ commutes with every $\rho_\nu(X^{(j)})$; it is a scalar, and $\rho_\nu(f)$ is a
multiple of $\Gamma$. By Fourier inversion, $f$ is then a multiple of $\sigma\mapsto\tr(\rho_\nu(\sigma^{-1})\Gamma)$. By a
theorem of Frobenius on the characters of the alternating group \cite[Thm.~2.5.13]{JK81}, this function vanishes off the
permutations of one cycle type, $(2n-1,2n-3,\dots,5,3)$, which for $n\ge3$ is not the type of an $N$-cycle. So $f$ vanishes
on the $N$-cycles as well as off them; it is zero (Lemma~\ref{lem:twocorner}).
\end{enumerate}

An expanded version of this paper, with further
results, is posted on Zenodo \cite{FV}; none of them is used here.

\section{Formanek's criterion in cycle form}\label{sec:criterion}

\subsection{Setting and standard facts}\label{subsec:setting}
As in Section~\ref{sec:intro}, $F$ is a field of characteristic zero, and $F\langle x_1,x_2,\dots\rangle$ is the free
\emph{non-unital} associative algebra on the letters $x_1,x_2,\dots$; its endomorphisms send each letter to a polynomial without
constant term. For $n\ge1$ fix a letter $x$ (for instance $x=x_1$), and let $\mathcal T:=(x^n)^T$ be the T-ideal generated by $x^n$; it does
not depend on the choice of $x$. By
Subsection~\ref{subsec:NH},
\[
d(n)=\min\{m: x_1\cdots x_m\in\mathcal T\}.
\]
Let $V_m$ be the space of multilinear polynomials in $x_1,\dots,x_m$ (its basis consists of the $m!$ words in which each of
$x_1,\dots,x_m$ occurs once) and $W_{n,m}:=V_m\cap\mathcal T$. We identify a partition $\lambda=(\lambda_1\ge\lambda_2\ge\dots)$
with its Young diagram, which has $\lambda_i$ cells in row $i$; the cell $(i,j)$ lies in row $i$ and column $j$, and its
\emph{content} is $j-i$. Exponents denote repeated parts, as in $(n^{n-1},n-1)=(n,\dots,n,n-1)$; we write $\ell(\lambda)$ for the number of rows, $\lambda'$ for the
conjugate (transposed) partition, and $\mu\subseteq\lambda$ for containment of diagrams. A \emph{removable corner} of $\lambda$
is a cell $a$ whose removal leaves a partition $\lambda-a$. We use the following standard facts; in (F2)--(F7), $N\ge1$ is
arbitrary, and from Subsection~\ref{subsec:criterion} on, $N=m+1$.

\begin{enumerate}[label=(F\arabic*),start=0,leftmargin=*]
\item \emph{Reduction to multilinear words} \cite[Prop.~1.8]{Ro22}. Since $\mathcal T$ is an ideal, \mbox{$d(n)\le m$} if and only if
$x_1\cdots x_m\in\mathcal T$ (if $x_1\cdots x_k\in\mathcal T$ with $k\le m$, then
\mbox{$x_1\cdots x_m=(x_1\cdots x_k)x_{k+1}\cdots x_m$} lies in $\mathcal T$). Since $\mathcal T$ is stable under renaming the
letters, which maps $x_1\cdots x_m$ onto every word of $V_m$, this holds if and only if \mbox{$W_{n,m}=V_m$}.

\item \emph{Spanning set.} For elements $y_1,\dots,y_n$ of $F\langle x_1,x_2,\dots\rangle$ let
\[
L(y_1,\dots,y_n):=\sum_{\pi\in S_n}y_{\pi(1)}\cdots y_{\pi(n)},
\]
and write $L$ for $L(y_1,\dots,y_n)$ when the $y_i$ are distinct letters. In characteristic $0$,
$\mathcal T=(L)^T$. Indeed, $L$ is the multilinear component of \mbox{$(y_1+\dots+y_n)^n$}, so $L\in\mathcal T$ (a T-ideal over an
infinite field contains the multihomogeneous components of its elements: replace each $y_i$ by $t_iy_i$ with scalars $t_i$,
and separate the components with a Vandermonde determinant); conversely, $x^n=L(x,\dots,x)/n!\in(L)^T$. The T-ideal $(L)^T$ is spanned by the elements
\mbox{$p\cdot L(h_1,\dots,h_n)\cdot q$}, where $p,q$ are words (possibly empty) and the $h_i$ have no constant term. Indeed, this
span lies in $(L)^T$, contains $L$, and is an ideal. It is closed under every endomorphism $\varphi$, because $\varphi$ maps
$p\cdot L(h_1,\dots,h_n)\cdot q$ to $\varphi(p)\cdot L(\varphi(h_1),\dots,\varphi(h_n))\cdot\varphi(q)$, the $\varphi(h_i)$ have
no constant term, and $\varphi(p)$, $\varphi(q)$ expand into words. By multilinearity the $h_i$ may be taken to be words. Each
such element is then multihomogeneous, so $W_{n,m}$ is spanned by the elements $p\cdot L(v_1,\dots,v_n)\cdot q$ in which
the $v_i$ are nonempty words, $p$ and $q$ are words (possibly empty), and each of $x_1,\dots,x_m$ occurs exactly once and no
other letter occurs. They have integer coefficients, so $\dim W_{n,m}$, and hence $d(n)$,
is the same for every field $F$ of characteristic zero.

\item \emph{Wedderburn decomposition} \cite[Thm.~4.12]{Ja78}, \cite[Ch.~1]{CST10}. For a partition $\lambda\vdash N$ let
$S^\lambda$ be the irreducible $F[S_N]$-module indexed by $\lambda$ (with $S^{(N)}$ the trivial module and $S^{(1^N)}$ the sign module), $\rho_\lambda\colon F[S_N]\to\End(S^\lambda)$ the representation and $\chi_\lambda$ its
character. The modules $S^\lambda$ are defined over $\mathbb{Q}$ and are absolutely irreducible. We have $F[S_N]=\bigoplus_{\lambda\vdash N}I_\lambda$
with $I_\lambda=F[S_N]e_\lambda$, where $e_\lambda=(\dim S^\lambda/N!)\sum_\sigma\chi_\lambda(\sigma)\sigma$ is a central primitive
idempotent, and $\rho_\lambda$ restricts to an isomorphism $I_\lambda\cong\End(S^\lambda)$; so each $I_\lambda$ is a simple algebra.
Fourier inversion reads
\[
f(\sigma)=\sum_\lambda\frac{\dim S^\lambda}{N!}\tr\bigl(\rho_\lambda(\sigma^{-1})\rho_\lambda(f)\bigr).
\]

\item \emph{Orthogonality.} Put $\langle f,f'\rangle:=\sum_\sigma f(\sigma)f'(\sigma)$. The permutations form an orthonormal
basis, so $\langle\,,\rangle$ is nondegenerate. The $I_\lambda$ are mutually orthogonal: $\langle f,f'\rangle$ is the coefficient
of $1$ in $f\cdot f'^\vee$, where $f'^\vee:=\sum_\sigma f'(\sigma)\sigma^{-1}$; the map $\vee$ is an anti-automorphism fixing each
$e_\lambda$, since the characters of $S_N$ are real, so $f\cdot f'^\vee\in I_\lambda I_\mu=0$ for $f\in I_\lambda$, $f'\in I_\mu$,
$\lambda\ne\mu$.

\item \emph{The ideal of the antisymmetrizer.} For $N\ge n+1$ regard $S_{n+1}\subseteq S_N$ (fixing $n+2,\dots,N$) and put
$E:=\sum_{\pi\in S_{n+1}}\sgn(\pi)\pi$. The two-sided ideal generated by $E$ is $I:=\bigoplus_{\ell(\lambda)\ge n+1}I_\lambda$.
Indeed, each $I_\lambda$ is simple, so this ideal is the sum of the $I_\lambda$ with $Ee_\lambda\ne0$, that is, with
$\rho_\lambda(E)\ne0$. As $E/(n+1)!$ is the idempotent of $F[S_{n+1}]$ for the sign character, $\rho_\lambda(E)\ne0$ if and only
if $\Res_{S_{n+1}}S^\lambda$ contains $\sgn$. By the branching rule (F5) below, applied $N-n-1$ times, this holds if and only if
$(1^{n+1})\subseteq\lambda$, that is, $\ell(\lambda)\ge n+1$.

\item \emph{Branching rule} \cite[Thm.~9.2(ii)]{Ja78}. $\Res_{S_{N-1}}S^\lambda\cong\bigoplus_a S^{\lambda-a}$, over the removable
corners $a$ of $\lambda$.

\item \emph{Burnside's theorem} \cite[Thm.~2.5(i)]{Et11}, which is part of the density theorem there. If $\rho\colon F[G]\to\End(U)$ is an absolutely irreducible representation
of a finite group $G$, then $\rho(F[G])=\End(U)$. (In \cite{Et11} the field is algebraically closed; over $F$ the two sides
have the same dimension as after extending scalars to an algebraic closure.)

\item \emph{Sign twist} \cite[Thm.~6.7]{Ja78}. The map $\varepsilon$, $\varepsilon(\sigma):=\sgn(\sigma)\sigma$, extends to an algebra automorphism of
$F[S_N]$ with $\varepsilon(I_\lambda)=I_{\lambda'}$, since $S^\lambda\otimes\sgn\cong S^{\lambda'}$. On $C_N$ it is the scalar
$(-1)^{N-1}$, because every $N$-cycle has sign $(-1)^{N-1}$.
\end{enumerate}

\subsection{Notation}
We compose permutations right to left, $(\sigma\tau)(i)=\sigma(\tau(i))$, and write $[k]:=\{1,\dots,k\}$ for $k\ge1$. The objects below
depend on $n$ and $N$; we write $B_n^{(N)}$ when the dependence on $N$ must be shown.
As in Subsection~\ref{subsec:outline}, $C_N$ is the span of the $N$-cycles, and $P(f):=\sum_\sigma f(\sigma)\sigma$, with $\sigma$
running over the $N$-cycles, is the coordinate projection onto $C_N$ (orthogonal for $\langle\,,\rangle$). Put
\[
A:=\bigoplus_{\ell(\lambda)\le n}I_\lambda,\qquad A':=\bigoplus_{\lambda_1\le n}I_\lambda,\qquad
B_n=B_n^{(N)}:=A\cap A'=\bigoplus_{\lambda\vdash N,\ \lambda\subseteq(n^n)}I_\lambda,
\]
the last sum being over the partitions of $N$ that fit in the $n\times n$ box. For a sequence $\alpha=(j_1,\dots,j_r)$ of
elements of $[m]$ we write $x_\alpha:=x_{j_1}\cdots x_{j_r}$ (the empty word if $r=0$). Orthogonal complements $\perp$ are taken
for $\langle\,,\rangle$.

\subsection{The criterion}\label{subsec:criterion}
\begin{proposition}[Formanek's criterion, cycle form]\label{prop:criterion}
Let $n\ge1$, $m\ge n$ and $N=m+1$. Then
\[
\dim W_{n,m}=m!-\dim\bigl(C_N\cap B_n^{(N)}\bigr).
\]
In particular, $d(n)\le m$ if and only if $C_N\cap B_n^{(N)}=0$.
\end{proposition}

\begin{proof}
Since $m\ge n$, we have $N\ge n+1$, so $S_{n+1}\subseteq S_N$ as in (F4). In steps (1)--(6) we show that the coordinate
projection $P$ maps the ideal $I$ of (F4) onto $\Phi(W_{n,m})$, where $\Phi$ is the isomorphism of step~(1); in step~(7) we pass
to orthogonal complements, using (F3) and the sign twist (F7).
\begin{enumerate}[label=(\arabic*),leftmargin=*]
\item \emph{Words and $N$-cycles.} Let $\Phi(x_{j_1}\cdots x_{j_m})$ be the $N$-cycle \mbox{$(N\,j_1\,\cdots\,j_m)$}, that is, $N\mapsto j_1\mapsto\dots\mapsto j_m\mapsto N$. This is a
bijection from the $m!$ multilinear words (the basis of $V_m$) onto the $N$-cycles, so it extends to a linear isomorphism
$\Phi\colon V_m\to C_N$. For example, for $m=3$,
$\Phi(x_2x_3x_1)$ is the $4$-cycle \mbox{$4\to2\to3\to1\to4$}.

\item \emph{Generators of $I$.} For $K\subseteq[N]$ with $|K|=n+1$, let $\Sym(K)\subseteq S_N$ be the group of the
permutations that fix every point outside $K$, and put \mbox{$E_K:=\sum_{\pi\in\Sym(K)}\sgn(\pi)\pi$}; thus $E=E_{[n+1]}$. Then $I$ is
spanned by the elements $\tau E_K$ with $\tau\in S_N$ and $|K|=n+1$. Indeed, $I$ is spanned by the elements $gEg'$ with
$g,g'\in S_N$, and $gEg'=(gg')(g'^{-1}Eg')=(gg')E_{g'^{-1}([n+1])}$; conversely, $\tau E_K=(\tau g)Eg^{-1}$ for any $g\in S_N$
with $g([n+1])=K$.

\item \emph{Reduction to $N$-cycles.} If no permutation $\tau\pi$ with $\pi\in\Sym(K)$ is an $N$-cycle, then $P(\tau E_K)=0$. Otherwise pick
$\pi_0\in\Sym(K)$ such that $\tau_0:=\tau\pi_0$ is an $N$-cycle; then $\tau E_K=\sgn(\pi_0)\tau_0E_K$, since $\tau=\tau_0\pi_0^{-1}$ and $\pi_0^{-1}E_K=\sgn(\pi_0)E_K$. So, to compute
$P(\tau E_K)$, we may assume that $\tau$ is an $N$-cycle.

\item \emph{Arcs.} For $w\in K$ let $\bar\tau(w):=\tau^k(w)$ for the least $k\ge1$ with $\tau^k(w)\in K$. Then $\bar\tau$ is an
$(n+1)$-cycle on $K$: the points of $K$ appear in some cyclic order along the single cycle $\tau$, and $\bar\tau$ sends each one
to the next. The \emph{arc} after $w\in K$ is $\Arc(w):=(\tau(w),\dots,\bar\tau(w))$. The $n+1$ arcs $\Arc(w)$, $w\in K$, are
nonempty and partition $[N]$. For $\pi\in\Sym(K)$, the permutation $\tau\pi$ agrees with $\tau$ off $K$ and sends each $w\in K$ to
$\tau(\pi(w))$; starting at $w$, it traverses $\Arc(\pi(w))$ and first returns to $K$ at $\bar\tau(\pi(w))$. Every point lies on an
arc and every arc ends in $K$, so every cycle of $\tau\pi$ passes through $K$, and the first-return map of $\tau\pi$ to $K$
(which sends $w\in K$ to the next point of $K$ on its $\tau\pi$-cycle) is $\bar\tau\pi$. The cycles of $\tau\pi$ correspond to
those of $\bar\tau\pi$, so $\tau\pi$ is an $N$-cycle if and only if $\bar\tau\pi$ is an $(n+1)$-cycle. For example, let $n=2$,
$N=5$, $\tau=\Phi(x_1x_2x_3x_4)=(5\,1\,2\,3\,4)$ and $K=\{1,3,4\}$. Along $\tau$, the points of $K$ come in the cyclic order
$1,3,4$, so $\bar\tau=(1\,3\,4)$; the arc after $1$ is $(2,3)$, the arc after $3$ is $(4)$, and the arc after $4$ is $(5,1)$.

\item \emph{Constant sign.} Let $\Pi:=\{\pi\in\Sym(K): \bar\tau\pi \text{ is an $(n+1)$-cycle}\}$. Since $\bar\tau$ is
itself an $(n+1)$-cycle, $\sgn(\bar\tau)=(-1)^n$, and for $\pi\in\Pi$,
\[
\sgn(\pi)=\sgn(\bar\tau)\sgn(\bar\tau\pi)=(-1)^n(-1)^n=+1.
\]
Therefore, for an $N$-cycle $\tau$, $P(\tau E_K)=\sum_{\pi\in\Pi}\tau\pi$.

\item \emph{The image of $I$.} The map $\pi\mapsto\theta:=\bar\tau\pi$ is a bijection from $\Pi$ onto the set of
$(n+1)$-cycles on $K$. By step~(4), the permutation $\tau\pi$, starting at $w\in K$, traverses $\Arc(\pi(w))$, which ends at
$\bar\tau(\pi(w))=\theta(w)$; repeating this, it traverses the arcs ending at $\theta(w),\theta^2(w),\dots$, in this order,
each from its first to its last point.
Write the arc that contains $N$ as $\Arc_0=(\zeta,N,\xi)$, where $\zeta$ is the part before $N$ and $\xi$ the part after $N$;
either may be empty. Number the other $n$ arcs $\Arc_1,\dots,\Arc_n$ in any fixed way. Read as a cycle starting at $N$, the
permutation $\tau\pi$ first finishes $\Arc_0$, which gives $\xi$, then traverses the arcs $\Arc_1,\dots,\Arc_n$ in some order, and
finally enters $\Arc_0$ again, which gives $\zeta$. (The word $\Phi^{-1}(\tau\pi)$ is read along the cycle starting just after
$N$; so the part $\xi$ after $N$ comes first in the word, and the part $\zeta$ before $N$ comes last.) Since each arc is determined by its last point, which lies in $K$, an $(n+1)$-cycle $\theta$ on $K$ amounts to a cyclic
order of the $n+1$ arcs,
that is, to an order of $\Arc_1,\dots,\Arc_n$ after $\Arc_0$; so as $\pi$ runs over $\Pi$, each of the $n!$ orders of
$\Arc_1,\dots,\Arc_n$ occurs exactly once. Hence
\[
\Phi^{-1}\Bigl(\sum_{\pi\in\Pi}\tau\pi\Bigr)=x_\xi\cdot L(x_{\Arc_1},\dots,x_{\Arc_n})\cdot x_\zeta ,
\]
where the $\Arc_i$ with $i\ge1$ are nonempty and do not contain $N$. By steps (2), (3) and (5), $P(I)$ is spanned by such sums, and each lies in
$\Phi(W_{n,m})$ by (F1); this proves the inclusion $P(I)\subseteq\Phi(W_{n,m})$. In the example of step~(4), \mbox{$\Arc_0=(5,1)$}, so $\zeta$ is empty and $\xi=(1)$, and the other two
arcs are $(2,3)$ and $(4)$. Moreover $\bar\tau\pi$ must be one of the two $3$-cycles $\bar\tau$, $\bar\tau^2$ on $K$, so
$\Pi=\{1,\bar\tau\}$; as \mbox{$\tau\bar\tau=(5\,1\,4\,2\,3)=\Phi(x_1x_4x_2x_3)$}, we get
\[
\Phi^{-1}P(\tau E_K)=x_1x_2x_3x_4+x_1x_4x_2x_3=x_1\cdot L(x_2x_3,x_4).
\]

\emph{Reverse inclusion.} Given $p\cdot L(v_1,\dots,v_n)\cdot q$ as in (F1), let $\tau:=\Phi(pv_1\cdots v_nq)$, an
\mbox{$N$-cycle}, and let $K$
consist of the indices of the last letters of the $v_i$ together with $z_0$, where $z_0$ is the index of the last letter of $p$,
or $z_0:=N$ if $p$ is empty. The $v_i$ are nonempty and all letters are distinct, so these are $n+1$ distinct points. Along
$\tau$ we meet $N,p,v_1,\dots,v_n,q,N$. So the arcs are (the index sequences of) $v_1,\dots,v_n$ and $\Arc_0=(q,N,p)$, or
$(q,N)$ if $p$ is empty; thus $\zeta=q$ and $\xi=p$. By step~(5) and the display above,
$\Phi(p\cdot L(v_1,\dots,v_n)\cdot q)=\sum_{\pi\in\Pi}\tau\pi=P(\tau E_K)\in P(I)$. Since these elements span $W_{n,m}$ by
(F1), $\Phi(W_{n,m})\subseteq P(I)$, and hence $P(I)=\Phi(W_{n,m})$.

\item \emph{Duality.} Let $f\in C_N$. Since $f(\sigma)=0$ unless $\sigma$ is an $N$-cycle, $\langle f,f'\rangle=\langle
f,P(f')\rangle$ for every $f'\in F[S_N]$. By step~(6), $f$ is orthogonal to $\Phi(W_{n,m})$ if and only if it is orthogonal to $I$.
Now $I^\perp=A$: indeed $A\subseteq I^\perp$ by (F3), and $\dim I^\perp=N!-\dim I$ because $\langle\,,\rangle$ is nondegenerate
(F3), while $N!-\dim I=\dim A$ because $F[S_N]=I\oplus A$ by (F2) and (F4). Thus $C_N\cap\Phi(W_{n,m})^\perp=C_N\cap A$. By (F7), $\varepsilon$ acts on $C_N$ as $\pm1$ and
maps $A$ onto $A'$ (since $\ell(\lambda)=\lambda'_1$); hence $C_N\cap A\subseteq A'$, and so $C_N\cap A=C_N\cap B_n$. The $N$-cycles form an orthonormal basis of
$C_N$, so $\langle\,,\rangle$ is nondegenerate on $C_N$, and $\dim C_N=(N-1)!=m!$. Since $\Phi(W_{n,m})\subseteq C_N$, the
space $C_N\cap\Phi(W_{n,m})^\perp$ is the orthogonal complement of $\Phi(W_{n,m})$ in $C_N$, so its dimension is
$m!-\dim\Phi(W_{n,m})$. Since $\Phi$ is an isomorphism,
\[
\dim W_{n,m}=m!-\dim\bigl(C_N\cap\Phi(W_{n,m})^\perp\bigr)=m!-\dim(C_N\cap B_n).
\]
\end{enumerate}
Finally, by (F0), $d(n)\le m$ if and only if $\dim W_{n,m}=m!$, that is, if and only if $C_N\cap B_n=0$.
\end{proof}

\begin{example}
Let $n=m=2$, so that $N=3$. By (F1), $W_{2,2}$ is spanned by \mbox{$L(x_1,x_2)=x_1x_2+x_2x_1$}, so $\dim W_{2,2}=1$ (and
$x_1x_2\notin\mathcal T$, that is, $d(2)>2$). On the group side, $C_3$ is spanned by the $3$-cycles $(1\,2\,3)$ and
$(1\,3\,2)$, which are $\Phi(x_1x_2)$ and $\Phi(x_2x_1)$ in step~(1), and $B_2^{(3)}=I_{(2,1)}$, since only $(2,1)$ fits
in the $2\times2$ box. Both $3$-cycles are even, so $f=a(1\,2\,3)+b(1\,3\,2)$ has the components \mbox{$(a+b)e_{(3)}$} and
\mbox{$(a+b)e_{(1^3)}$} in the other two blocks; hence $C_3\cap I_{(2,1)}$ is spanned by $(1\,2\,3)-(1\,3\,2)$, and
\mbox{$m!-\dim(C_3\cap B_2^{(3)})=2-1=\dim W_{2,2}$}, as Proposition~\ref{prop:criterion} says. In step~(7), this is the
statement that the orthogonal complement of $\Phi(W_{2,2})=F\cdot\bigl((1\,2\,3)+(1\,3\,2)\bigr)$ in $C_3$ is
$F\cdot\bigl((1\,2\,3)-(1\,3\,2)\bigr)$.
\end{example}

\begin{remark}[the direction used]\label{rem:forward}
Theorems~\ref{thm:A} and~\ref{thm:B} use only the direction in which $C_N\cap B_n=0$ implies $d(n)\le m$. This direction needs only
step~(7) and the inclusion $P(I)\subseteq\Phi(W_{n,m})$ of step~(6); the reverse inclusion gives the exact formula.
\end{remark}

\begin{remark}[Razmyslov's bound]\label{rem:razmyslov}
Take $m=n^2$, so that $N=n^2+1$. No partition of $N$ fits in the $n\times n$ box, so $B_n^{(N)}=0$, and
Proposition~\ref{prop:criterion} gives $d(n)\le n^2$.
\end{remark}

\section{Proof of Theorem~\ref{thm:A}}\label{sec:A}

\begin{lemma}[the block of the square]\label{lem:onecorner}
For $n\ge2$, $C_{n^2}\cap B_n^{(n^2)}=0$.
\end{lemma}

\begin{proof}
Let $N=n^2$. The only partition of $n^2$ that fits in the $n\times n$ box is the square $\beta:=(n^n)$, so $B_n^{(N)}=I_\beta$. Let $f\in C_N\cap
I_\beta$ and $M:=\rho_\beta(f)$. Then $f(\sigma)=(\dim S^\beta/N!)\tr(\rho_\beta(\sigma^{-1})M)$ by (F2). A permutation $\sigma$
in the stabilizer $S_{N-1}$ of the point $N$ is not an $N$-cycle, so $f(\sigma)=0$. The square has only one removable corner, so
by (F5) $\Res_{S_{N-1}}S^\beta$ is the module $S^\nu$ with $\nu=(n^{n-1},n-1)$, which is absolutely irreducible by (F2); by
(F6), the $\rho_\beta(\sigma)$ with $\sigma\in S_{N-1}$ span $\End(S^\beta)$. Hence $\tr(ZM)=0$ for every $Z\in\End(S^\beta)$, so $M=0$ and $f=0$.
\end{proof}

\begin{proof}[Proof of Theorem~\ref{thm:A}]
Lemma~\ref{lem:onecorner} and Proposition~\ref{prop:criterion} with $m=n^2-1\ge n$ give $d(n)\le n^2-1$.
\end{proof}

\section{Proof of Theorem~\ref{thm:B}}\label{sec:B}

We follow the outline in Subsection~\ref{subsec:outline}(ii).

\subsection{Reduction to one block}\label{subsec:reduction}
Let $N=n^2-1$. A partition of $N$ that fits in the $n\times n$ box leaves out exactly one cell of the box, and this cell must
be $(n,n)$, the only cell of the box whose removal leaves a diagram. Thus the only partition of $N$ that fits is $\nu:=(n^{n-1},n-1)$, the box with the cell $(n,n)$ removed. Hence $B_n=I_\nu$,
and by Proposition~\ref{prop:criterion} it suffices to show $C_N\cap I_\nu=0$. By (F1), $d(n)$ does not depend on $F$, so
from now on we take $F=\mathbb{C}$, which \ref{T3} below requires. For example, for $n=3$, $\nu=(3,3,2)$ and $N=8$; the contents of its cells are $0,1,2$ in the first row, $-1,0,1$ in
the second and $-2,-1$ in the third, and its two removable corners are the cells $(2,3)$ and $(3,2)$, of contents $+1$ and
$-1$. For $f\in I_\nu$ put $M:=\rho_\nu(f)$; then, by (F2),
\[
f(\sigma)=\frac{\dim S^\nu}{N!}\tr\bigl(\rho_\nu(\sigma^{-1})M\bigr).\tag{$*$}
\]

\subsection{Facts from the representation theory of \texorpdfstring{$S_N$}{S\_N}}\label{subsec:facts}
Besides (F2) and (F5)--(F7), we use three facts. Below, $\rho$ is $\rho_\lambda$ for a partition $\lambda\vdash N$ (from
Subsection~\ref{subsec:stars} on, $\lambda=\nu$). A standard tableau $T$ of shape $\lambda$ is a filling of the diagram of
$\lambda$ with $1,\dots,N$, increasing along the rows and down the columns; for an entry $k$, $c_T(k)$ is the content of the
cell of $k$ in $T$, and $T_{\le k}$ is the shape formed by the cells of $1,\dots,k$. The \emph{Young basis} $\{v_T\}$ of
$S^\lambda$ is defined, up to nonzero scalars, by following $T$ down the chain $S_N\supset S_{N-1}\supset\dots\supset S_1$
\cite{OV96}, \cite[\S3.2]{CST10}: by the branching rule (F5), $S^\lambda$ contains exactly one $F[S_{N-1}]$-submodule
isomorphic to $S^{T_{\le N-1}}$, this submodule contains exactly one $F[S_{N-2}]$-submodule isomorphic to $S^{T_{\le N-2}}$,
and so on; the last of these submodules is a line, and $v_T$ spans it. (This is not Young's natural basis of polytabloids
\cite{Ja78}.)

\begin{enumerate}[label=(T\arabic*),leftmargin=*]
\item\label{T1} \emph{Jucys--Murphy elements.} Put $J_k:=\sum_{i<k}(i\,k)$. In the Young (Gelfand--Tsetlin) basis $\{v_T\}$,
$\rho(J_k)v_T=c_T(k)v_T$ \cite{Ju74,Mu81,OV96}; so $\rho(J_k)$ is diagonal in this basis. The basis is adapted to restriction:
for each removable corner $a$ of $\lambda$, the $v_T$ with $N$ in the cell $a$ span the summand $S^{\lambda-a}$ of
$\Res_{S_{N-1}}S^\lambda$, and $J_N$ acts on it by the content of $a$.

\item\label{T2} \emph{Young's seminormal form} \cite{OV96}, \cite[Ch.~3]{CST10}. Let $s:=(k\,k{+}1)$, and let $sT$ be $T$ with
the entries $k$ and $k+1$ exchanged. In the basis $\{v_T\}$ of \ref{T1}, $\rho(s)v_T\in\operatorname{span}(v_T,v_{sT})$, with the $v_{sT}$ term only if $sT$
is standard. Moreover $sJ_{k+1}-J_ks=1$ (since $s\,(i\,k{+}1)=(i\,k)\,s$ for $i<k$, and $s\,(k\,k{+}1)=1$). Write
$\rho(s)_{TT}$ for the $v_T$-coefficient of $\rho(s)v_T$. Taking the $v_T$-coefficient of this identity applied to $v_T$
($\rho(J_k)$ and $\rho(J_{k+1})$ are diagonal by \ref{T1}) gives $(c_T(k+1)-c_T(k))\,\rho(s)_{TT}=1$, so
\[
\rho(s)_{TT}=\frac{1}{c_T(k+1)-c_T(k)}.
\]

\item\label{T3} \emph{The operator $\Gamma$ and Frobenius' theorem.} Let $\lambda=\lambda'$ with $N\ge2$. Then, by (F7), $S^\lambda\otimes\sgn\cong
S^\lambda$, so there is an invertible operator $\Gamma$ on $S^\lambda$ with $\Gamma\rho(\sigma)=\sgn(\sigma)\rho(\sigma)\Gamma$
for all $\sigma\in S_N$. The operator $\Gamma^2$ commutes with $\rho(S_N)$, so by Schur's lemma it is a nonzero scalar, and
over $\mathbb{C}$ we rescale $\Gamma$ so that $\Gamma^2=\mathrm{id}$. The operator $\Gamma$ commutes with $\rho(A_N)$, where $A_N$ is the alternating group, and the
restriction of $S^\lambda$ to $A_N$ is the sum $U_+\oplus U_-$ of two non-isomorphic irreducible modules, with characters
$\chi_\pm$ \cite[\S2.5]{JK81}. By Schur's lemma, $\Gamma$ maps each $U_\pm$ into itself and acts on it by a scalar, which is $1$ or $-1$ since
$\Gamma^2=\mathrm{id}$. The two scalars differ, because $\Gamma\ne\pm\mathrm{id}$ (it anticommutes with $\rho((1\,2))\ne0$). So
the function
\[
\psi(\sigma):=\tr(\rho(\sigma)\Gamma)
\]
vanishes on odd permutations (there $\Gamma\rho(\sigma)=-\rho(\sigma)\Gamma$), and $\psi=\pm(\chi_+-\chi_-)$ on $A_N$. Let
$(1,1),\dots,(r,r)$ be the diagonal cells of $\lambda$. The \emph{hook} of a diagonal cell $(c,c)$ consists of the cells
$(c,j)$ of $\lambda$ with $j\ge c$ and $(i,c)$ of $\lambda$ with $i>c$; as $\lambda=\lambda'$, the hook of $(i,i)$ has
$h_{ii}=2(\lambda_i-i)+1$ cells. These numbers are odd, and they are distinct because $\lambda_i-i$ decreases strictly with
$i$. Their sum is $N$, since every cell $(a,b)$ of $\lambda$ lies in the hook of exactly one diagonal cell, namely $(c,c)$
with $c=\min(a,b)$ (a cell of $\lambda$, since $c\le a$ and $c\le b$). So $h(\lambda):=(h_{11},\dots,h_{rr})$ is a cycle type of $S_N$. \emph{Frobenius' theorem} \cite{Fr01}, \cite[Thm.~2.5.13]{JK81}
says that $\chi_+(\sigma)=\chi_-(\sigma)$ unless $\sigma$ has cycle type $h(\lambda)$. So $\psi$ vanishes on every permutation
whose cycle type is not $h(\lambda)$.
\end{enumerate}

\begin{example}
Let $\lambda=(2,1)$, so that $N=3$, and let $T_{12}$ and $T_{13}$ be the standard tableaux with first rows $1,2$ and $1,3$. In the
basis $v_{T_{12}},v_{T_{13}}$, \ref{T1} gives \mbox{$\rho((1\,2))=\rho(J_2)=\operatorname{diag}(1,-1)$}, and \ref{T2} gives the diagonal
entries $1/(-1-1)=-\frac12$ and $1/(1+1)=\frac12$ of $\rho((2\,3))$. With the $v_T$ scaled as in Young's orthogonal form
\cite[\S3.4.2]{CST10},
\[
\rho((1\,2))=\begin{pmatrix}1&0\\0&-1\end{pmatrix},\qquad
\rho((2\,3))=\begin{pmatrix}-\frac12&\frac{\sqrt3}2\\[2pt]\frac{\sqrt3}2&\frac12\end{pmatrix}.
\]
An operator that anticommutes with both matrices has zero diagonal entries (by the first) and opposite off-diagonal entries
(by the second), so it is a multiple of $\Gamma_0:=\left(\begin{smallmatrix}0&1\\-1&0\end{smallmatrix}\right)$. As
$\Gamma_0^2=-\mathrm{id}$, the operator $\Gamma$ of \ref{T3} is $\pm i\,\Gamma_0$; this is why \ref{T3} needs $\mathbb C$.
For $\Gamma=i\,\Gamma_0$, $\psi$ vanishes on the identity and on the transpositions, and $\psi((1\,2\,3))=-i\sqrt3=
-\psi((1\,3\,2))$, as Frobenius' theorem says for $h(\lambda)=(3)$.
\end{example}

\subsection{Star elements and permutations with a fixed point}\label{subsec:stars}
For $j\in[N]$ define the \emph{star element} $X^{(j)}:=\sum_{i\ne j}(i\,j)\in F[S_N]$, the sum of the $N-1$ transpositions that
move $j$, and put $Y_j:=\rho(X^{(j)})$. Note that $X^{(N)}=J_N$ and $gX^{(j)}g^{-1}=X^{(g(j))}$ for $g\in S_N$. The shape $\nu$
has exactly two removable corners, $a_+:=(n-1,n)$ of content $+1$ and $a_-:=(n,n-1)$ of content $-1$; the other row ends
$(i,n)$, $i<n-1$, are not removable, since the next row has the same length. Write $\mathrm{Stab}(j)\subseteq S_N$ for the
stabilizer of the point $j$, so $\Stab(N)=S_{N-1}$. By (F5), \mbox{$\Res_{\Stab(N)}S^\nu=S^{\nu-a_+}\oplus S^{\nu-a_-}$}, and by \ref{T1}
$Y_N=\rho(J_N)$ acts on the two summands by $+1$ and $-1$; let $Q_\pm:=(\mathrm{id}\pm Y_N)/2$ be the projections onto them.

\begin{lemma}[vanishing on permutations with a fixed point]\label{lem:anticommute}
Let $f\in I_\nu$ and $M:=\rho(f)$. Then $f$ vanishes on every permutation with a fixed point if and only if
$MY_j=-Y_jM$ for all $j\in[N]$.
\end{lemma}

\begin{proof}
Every permutation with a fixed point lies in some $\Stab(j)$. We show that $f$ vanishes on $\Stab(j)$ if and only if
$MY_j=-Y_jM$, first for $j=N$. The two summands $S^{\nu-a_\pm}$ are absolutely irreducible and non-isomorphic, since \mbox{$\nu-a_+\ne\nu-a_-$}. By (F2) for
$\Stab(N)$, the
isotypic blocks of $F[\Stab(N)]$ of these two shapes map onto $\End(S^{\nu-a_+})$ and $\End(S^{\nu-a_-})$ respectively, and each
acts as zero on the other summand; so \mbox{$\rho(F[\Stab(N)])=\End(S^{\nu-a_+})\oplus\End(S^{\nu-a_-})$}, the operators that preserve
both summands. By ($*$), $f$ vanishes on $\Stab(N)$ if and only if $\tr(ZM)=0$ for every $Z\in\rho(F[\Stab(N)])$. These $Z$
are the operators $Q_+Z_1Q_++Q_-Z_2Q_-$ with $Z_1,Z_2\in\End(S^\nu)$, and \mbox{$\tr(Q_\pm Z_iQ_\pm M)=\tr(Z_i\,Q_\pm MQ_\pm)$}. Since
the trace form $(Z,Z')\mapsto\tr(ZZ')$ on $\End(S^\nu)$ is nondegenerate, the condition holds if and only if
$Q_+MQ_+=0=Q_-MQ_-$. Since $Y_N$ is an
involution and $Q_\pm$ are the projections onto its $\pm1$-eigenspaces, this says that $M$ maps each eigenspace of $Y_N$ into the
other, that is, $MY_N=-Y_NM$.

\emph{Transfer to $\Stab(j)$.} Take $g\in S_N$ with $g(N)=j$, and apply the result for $\Stab(N)$ to $g^{-1}fg\in I_\nu$: its coefficient
function is $\sigma\mapsto f(g\sigma g^{-1})$, and $\rho(g^{-1}fg)=\rho(g)^{-1}M\rho(g)$. Since $g\Stab(N)g^{-1}=\Stab(j)$, $f$ vanishes on $\Stab(j)$ if and only if $\rho(g)^{-1}M\rho(g)$ anticommutes with $Y_N$,
that is, if and only if $M$ anticommutes with $\rho(g)Y_N\rho(g)^{-1}=\rho(gJ_Ng^{-1})=Y_j$.
\end{proof}

\begin{lemma}[the commutant of the star elements]\label{lem:claimS}
For $n\ge2$, the only operators on $S^\nu$ that commute with all $Y_j$ are the scalars.
\end{lemma}

\begin{proof}
Let $\mathcal{Y}$ be the unital algebra generated by $Y_1,\dots,Y_N$. Since $\rho(g)Y_j\rho(g)^{-1}=Y_{g(j)}$, conjugation by
$\rho(g)$ maps $\mathcal{Y}$ onto itself. We show that $\rho(s)\in\mathcal Y$ for \mbox{$s:=(N{-}1\;N)$}.
\begin{enumerate}[label=(\alph*),leftmargin=*]
\item $Y_N=\rho(J_N)\in\mathcal{Y}$, so $Q_\pm\in\mathcal{Y}$.
\item $Y_{N-1}=\rho(J_{N-1})+\rho(s)$, since $X^{(N-1)}=\sum_{i<N-1}(i\;N{-}1)+(N{-}1\;N)$.
\item $D:=Q_+Y_{N-1}Q_++Q_-Y_{N-1}Q_-\in\mathcal{Y}$ is the block-diagonal part of $Y_{N-1}$; by \ref{T1}, $Q_\pm$ are the
coordinate projections onto the span of the $v_T$ with $N$ in $a_\pm$. By \ref{T1}, $\rho(J_{N-1})$ is diagonal, so it is its own
block-diagonal part. By \ref{T2}, the only off-diagonal entries of $\rho(s)$ that can be nonzero are the $v_{sT}$-coefficients of $\rho(s)v_T$ with $T$
and $sT$ both standard. Then the cells of $N$ in $T$ and in $sT$ are both removable corners of $\nu$ (in a standard tableau, the cell of the
largest entry is always a removable corner), and they are different, so
$N-1$ and $N$ occupy $a_+$ and $a_-$ in some order, and $c_T(N)\ne c_{sT}(N)$. So each nonzero off-diagonal entry of $\rho(s)$
links the two eigenspaces of $Y_N$, and the block-diagonal part of $\rho(s)$ is its diagonal. Hence, by \ref{T2} with $k=N-1$,
\[
Dv_T=\Bigl(c_T(N-1)+\frac{1}{c_T(N)-c_T(N-1)}\Bigr)v_T .
\]
\item Let $u=\pm1$ be the content of the cell of $N$ in $T$ (this cell is $a_+$ or $a_-$). Then $N-1$ is in a removable
corner of $\nu$ minus that cell: for $N$ in $a_+$, the corners of \mbox{$\nu-a_+=(n^{n-2},n-1,n-1)$} are $a_-$ (content $-1$) and, if
$n\ge3$, $(n-2,n)$ (content $+2$); for $N$ in $a_-$, those of \mbox{$\nu-a_-=(n^{n-1},n-2)$} are $a_+$ (content $+1$) and, if $n\ge3$,
$(n,n-2)$ (content $-2$). So $c_T(N-1)$ is $-u$ or $2u$, and in both cases
$c_T(N-1)\bigl(c_T(N)-c_T(N-1)\bigr)=-2u^2=-2$. Hence the coefficient in (c) is
\[
c_T(N-1)+\frac{1}{c_T(N)-c_T(N-1)}=c_T(N-1)-\frac{c_T(N-1)}{2}=\frac12\,c_T(N-1),
\]
so $Dv_T=\tfrac12c_T(N-1)v_T$ for every $T$, that is, $D=\tfrac12\rho(J_{N-1})$. So $\rho(J_{N-1})=2D\in\mathcal{Y}$, and
\mbox{$\rho(s)=Y_{N-1}-2D\in\mathcal{Y}$}. For example, for $n=3$ (so $\nu=(3,3,2)$ and $N=8$), if $8$ is in $a_+=(2,3)$, then $7$
is in $(1,3)$ or in $a_-=(3,2)$, of contents $2$ and $-1$, and $c_T(7)\bigl(c_T(8)-c_T(7)\bigr)$ is $2\cdot(1-2)=-2$ or
$(-1)\cdot(1+1)=-2$.
\item By (d), $\rho((N{-}1\;N))\in\mathcal{Y}$; conjugating, $\rho(t)\in\mathcal{Y}$ for every transposition $t$. So
$\mathcal{Y}\supseteq\rho(F[S_N])=\End(S^\nu)$ by (F6), and the commutant of $\mathcal Y$ is $\mathbb{C}\cdot\mathrm{id}$.\qedhere
\end{enumerate}
\end{proof}

\subsection{Conclusion of the proof}\label{subsec:end}
\begin{lemma}[the block of $\nu$]\label{lem:twocorner}
For every $n\ge3$, $C_{n^2-1}\cap B_n^{(n^2-1)}=0$.
\end{lemma}

\begin{proof}
Let $f\in C_N\cap I_\nu$, with $M=\rho(f)$ as in ($*$). A permutation with a fixed point is not an $N$-cycle, so $M$ anticommutes with
every $Y_j$ by Lemma~\ref{lem:anticommute}. The shape $\nu$ equals its conjugate (columns $1,\dots,n-1$ have length $n$ and
column $n$ has length $n-1$), so the operator $\Gamma$ of \ref{T3} exists. Each $X^{(j)}$ is a sum of transpositions, so
$\Gamma$ anticommutes with every $Y_j$. Then $\Gamma M$ commutes with every $Y_j$, so $\Gamma M=\kappa\,\mathrm{id}$ for some
$\kappa\in\mathbb C$ by Lemma~\ref{lem:claimS}, and $M=\kappa\Gamma$, as $\Gamma^2=\mathrm{id}$. By ($*$), $f(\sigma)=\kappa'\psi(\sigma^{-1})$ with $\kappa':=\kappa\dim S^\nu/N!$, and
\ref{T3} shows that $f$ vanishes on every permutation whose cycle type is not $h(\nu)$ (as $\sigma^{-1}$ has the cycle type of
$\sigma$). The diagonal cells of $\nu$ are $(i,i)$
for $i=1,\dots,n-1$, with $h_{ii}=2(n-i)+1$, so
\[
h(\nu)=(2n-1,2n-3,\dots,5,3);
\]
for example, $h(\nu)=(5,3)$ for $n=3$. For $n\ge3$ this type has $n-1\ge2$ cycles, so it is not the type of an $N$-cycle; hence
$f$ vanishes on the $N$-cycles. Since $f\in C_N$, $f$ also vanishes off the $N$-cycles. So $f=0$.
\end{proof}

\begin{proof}[Proof of Theorem~\ref{thm:B}]
Lemma~\ref{lem:twocorner} and Proposition~\ref{prop:criterion} with $m=n^2-2\ge n$ give $d(n)\le n^2-2$ over $\mathbb C$,
and hence over every field $F$ of characteristic zero by (F1).
\end{proof}

\begin{remark}[the case $n=2$]
For $n=2$, $\nu=(2,1)$ and $h(\nu)=(3)$ is the cycle type of the $3$-cycles, so the last step of the proof of
Lemma~\ref{lem:twocorner} fails, and the lemma itself is false for $n=2$. Indeed
$f=(1\,2\,3)-(1\,3\,2)$ is a nonzero element of $C_3\cap I_{(2,1)}$ (see the example after Proposition~\ref{prop:criterion}),
so $d(2)>2$ by Proposition~\ref{prop:criterion}; with Theorem~\ref{thm:A}, $d(2)=3$.
\end{remark}

\section*{Acknowledgments}
The author thanks Claudio Procesi for helpful suggestions.
This work was done overwhelmingly by AI models and supervised by the author. The models were Anthropic's Claude and OpenAI's ChatGPT. The author takes full responsibility for the paper.

\end{document}